\documentclass[10pt]{amsart}
\usepackage{amssymb,latexsym,amsmath,amsfonts}
\usepackage[mathscr]{eucal}

\numberwithin{equation}{section}
\theoremstyle{definition}
\newtheorem{definition}{Definition}[section]

\theoremstyle{plain}
\newtheorem{theorem}[definition]{Theorem}
\newtheorem{lemma}[definition]{Lemma}

\newcommand{\eps}{\varepsilon}
\newcommand{\sma}{\epsilon}

\newcommand{\zt}{\zeta}
\newcommand{\zbar}{\overline{z}}
\newcommand{\tht}{\theta}
  
\newcommand{\unal}{\mathscr{A}}
\newcommand{\nat}{\mathbb{N}}

\newcommand{\Leb}{\mathfrak{m}}
\newcommand{\exep}{\mathcal{E}}
\newcommand{\zee}{\mathfrak{z}}
\newcommand{\wee}{\mathfrak{w}}
\newcommand{\fee}{\mathfrak{F}}
\newcommand{\pow}{\varDelta}

\newcommand\partl[2]{\frac{\partial{#1}}{\partial{#2}}}

\newcommand\ball[2]{\overline{D(}{#1};{#2})}
\newcommand\sect[2]{\overline{S(}{#1};{#2})}
\newcommand{\OM}{\Omega}
\newcommand{\Dsc}{\overline{D}}
\newcommand\cyl[1]{\{(z,w):|z|\leq #1\}}

\newcommand{\smoo}{\mathcal{C}}
\newcommand{\hol}{\mathcal{O}}
\newcommand{\poly}{\mathscr{P}}

\newcommand{\er}{{\rm Re}}
\newcommand{\mi}{{\rm Im}}

\newcommand{\mideal}{\mathscr{M}[\mathscr{P}(K;\mathbb{C}^n)]}

\newcommand{\CC}{\mathbb{C}^2}

\newcommand{\cplx}{\mathbb{C}}

\newcommand{\srf}{\mathcal{S}}   

\newcommand{\presrf}{\mathfrak{S}}

\begin{document}

\title[Local Polynomial convexity]{Polynomial approximation, local polynomial \\
convexity, and degenerate CR singularities}
\author{Gautam Bharali}
\thanks{\bf To appear in {\em Journal of Functional Analysis}}
\address{Department of Mathematics, University of Michigan, MI 48109}
\email{bharali@umich.edu}
\keywords{Complex tangency, CR singularity, polynomial approximation, polynomially convex}
\subjclass[2000]{Primary 30E10, 32E20, 46J10}

\begin{abstract} 
We begin with the following question: given a closed disc $\Dsc\Subset\cplx$ and a complex-valued
function $F\in\smoo(\Dsc)$, is the uniform algebra on $\Dsc$ generated by $z$ and $F$ equal
to $\smoo(\Dsc)$ ? When $F\in\smoo^1(D)$,
this question is complicated by the presence of points in the surface
$\srf:={\rm graph}_{\Dsc}(F)$ that have complex tangents. Such points are called
CR singularities. Let $p\in\srf$ be a CR singularity at which the order of contact of the
tangent plane with $\srf$ is greater than $2$; i.e. a degenerate CR singularity. We provide sufficient
conditions for $\srf$ to be locally polynomially convex at the degenerate singularity $p$. This is useful
because it is essential to know whether $\srf$ is locally polynomially convex at a CR singularity in
order to answer the initial question. To this end, we also present a general theorem on the uniform
algebra generated by $z$ and $F$, which we use in our investigations. This result may be of
independent  interest because it is applicable even to non-smooth, complex-valued $F$.
\end{abstract}
\maketitle

\section{Introduction and statement of results}\label{S:intro}

One of the concerns of this paper is to study the following question: given a closed disc 
$\Dsc\Subset\cplx$ and a complex-valued function $F\in\smoo(\Dsc)$, when is the uniform algebra
on $\Dsc$ generated by $z$ and $F$ equal to $\smoo(\Dsc)$ ? A necessary condition for a positive
answer to this question is that ${\rm graph}_{\Dsc}(F)\subset\CC$ must be polynomially convex.
A compact subset $K\subset\cplx^n$ is said to be {\bf polynomially convex} if for each point
$\zt\notin K$, there exists a holomorphic polynomial $P$ such that $P(\zt)=1$ and $\sup_K|P|<1$.
The compact $K$ is said to be {\bf locally polynomially convex} at a point $p\in K$ if there
exists a closed ball $\mathbb{B}(p)$ centered at $p$ such that $K\cap\mathbb{B}(p)$
is polynomially convex. In general, it is difficult to determine whether a given compact
$K\subset\cplx^n$ is polynomially convex when $n>1$, but questions of polynomial convexity 
arise repeatedly in connection with function theory. There is a considerable
body of work concerning the (local) polynomial convexity of smooth surfaces in $\cplx^n$. The
references associated with the smooth case are too numerous to list here; instead, the reader is
referred to the survey \cite{sGarcia:phsd00}. In the instances discussed in that survey, one also
obtains positive answers to the question presented above. In contrast, very little is known when
$F$ is non-smooth --- either about the polynomial convexity of ${\rm graph}_{\Dsc}(F)$, or about
the question asked above --- beyond Mergelyan's result \cite{mergelyan:uafcv52}.
Mergelyan's result, however, is only applicable when $F$ is real-valued. Part of the intention
of this paper is to provide a sufficient condition on a {\em complex-valued} $F\in\smoo(\Dsc)$
for the question posed above to have an affirmative answer. This result is Theorem
\ref{T:uaprx} stated below. We note that our sufficient condition is stated in terms of the
value-distribution of $F$, which is easy to understand, and may be applied to concrete
situations. One such application is Theorem \ref{T:pcvx} below. 
\smallskip

Before stating Theorem \ref{T:uaprx}; we need to introduce some notation. In what
follows, $\Dsc$ will denote any closed disc in $\cplx$, while $\ball{a}{r}$ will denote 
the closed disc of radius $r$ centered at $a\in\cplx$. Unless explicitly stated otherwise, all
continuous functions will be assumed to be {\em complex-valued}. If $K$ is a compact subset of 
$\cplx$, and $\phi_1,\phi_2,\dots,\phi_N$ are continuous functions on $K$, 
$[\phi_1,\phi_2,\dots,\phi_N]_K$ is defined as
\begin{multline}
[\phi_1,\phi_2,\dots,\phi_N]_K \ := \ 
\{f\in\smoo(K):f \ \text{can be approximated uniformly on $K$}\notag\\
		\text{by complex polynomials in $\phi_1,\phi_2,\dots,\phi_N$}\}.\notag
\end{multline}
An {\bf open sector with vertex at $\boldsymbol{a}$}, denoted by $S(a;I)$, is the set
\[
S(a;I) \ := \ \{a+re^{i\theta}:r>0, \text{and} \ \theta\in I\},
\]
where $I$ is an open subinterval of $[-2\pi,2\pi)$ with ${\rm length}(I)<2\pi$, and $a\in\cplx$.
Having established our notation, we can now state our first result.
\smallskip

\begin{theorem}\label{T:uaprx} Let $F$ be a complex-valued continuous function on a closed disc 
$\Dsc\Subset\cplx$. Suppose that there is a set $E\subset\Dsc$ having zero Lebesgue measure such
that $F^{-1}\{F(\zt)\}$ is at most countable $\forall\zt\in\Dsc\setminus E$. Furthermore, suppose
that for each $\zt\in\Dsc\setminus E$, there exists an open sector $S(0;I_\zt)$ with vertex at 
$0\in\cplx$ such that
\begin{equation}\label{E:sector}
(z-\zt)\{F(z)-F(\zt)\} \ \in \ S(0;I_\zt) \quad\forall z\in\Dsc\setminus F^{-1}\{F(\zt)\}.
\end{equation}
Then, $[z,F]_{\Dsc}=\smoo(\Dsc)$.
\end{theorem}

The reader will notice that Weierstrass's approximation theorem is a special case of the above
theorem: when $F(z)=\zbar$, \eqref{E:sector} is satisfied by taking $S(0;I_\zt)$, for each
$\zt\in\Dsc$, to be a {\em fixed} sector containing the positive real axis.
We remark here that the proof of the above theorem is reminiscent of the early work of 
Wermer --- see, for instance, \cite[Theorem 1]{wermer:ad64} --- on questions of the sort 
considered in this paper. The crucial difference between those results and Theorem \ref{T:uaprx}
is that the sectors $S(0;I_\zt)$ occurring herein are allowed to have interior angles that are
{\em greater than $\pi$}. In results such as \cite[Theorem 1]{wermer:ad64}, the conclusion
$[z,F]_{\Dsc}=\smoo(\Dsc)$ is obtained under the assumption that $F$ is injective on $\Dsc$, and
the methods used in those results work only if the quantities occurring in \eqref{E:sector} lie 
in small subsets of a half-plane. It is for these reasons that Theorem \ref{T:uaprx} is a more 
general result.
\smallskip

Another objective of this paper is to further investigate the smooth case. Let $\srf$ be a smooth real
surface $\srf$ in $\cplx^n, \ n>1$. A point $p\in\srf$ is said to be {\bf totally real} if the tangent
plane $T_p(\srf)$ at $p$ is not a complex line. A point on $\srf$ that is not totally real will be 
called a {\bf CR singularity}. At a totally real point $p\in\srf$, the surface $\srf$ is locally 
polynomially convex. Contrast this with a CR singularity $p\in\srf\subset\CC$ when the order of 
contact of $T_p(\srf)$ with $\srf$ equals $2$. Since $T_p(\srf)$ is a complex tangent, there exist
holomorphic coordinates $(z,w)$ centered at $p$ such that $\srf$
is locally given by an equation of the form $w=|z|^2+\gamma(z^2+\zbar^2)+G(z)$, where $\gamma\geq 0$,
$G(z)=O(|z|^3)$, and three distinct situations arise. In Bishop's terminology, the CR singularity
$p=(0,0)$ is said to be elliptic if $0\leq\gamma<1/2$, parabolic if $\gamma=1/2$, and hyperbolic if
$\gamma>1/2$. Bishop showed \cite{bishop:dmcEs65}, among other things, that if $p$ is elliptic,
then $\srf$ is {\em not} locally polynomially convex. Much later, Forstneri{\v{c}} \& Stout
\cite{forstnericStout:ncpcs91} showed (also refer to \cite{stout:pcnhp86} by Stout) that if 
$p\in\srf$ is an isolated, hyperbolic CR singularity, then $\srf$ is locally polynomially convex
at $p$. Furthermore, in the hyperbolic case, writing $F(z)=|z|^2+\gamma(z^2+\zbar^2)+G(z)$, it
has been shown in \cite{forstnericStout:ncpcs91} that, for a small $\eps>0$, 
$[z,F]_{\{|z|\leq\eps\}}=\smoo(\ball{0}{\eps})$. The more involved case $\gamma=1/2$ has been
studied in \cite{joricke:lphdnipp97}.
\smallskip

This raises the question: what can be said about the polynomial convexity of a surface $\srf$ if
the order of contact of $T_p(\srf)$ with $\srf$ at a CR singularity $p$ is {\em greater} than $2$ ?
We will call such a CR singularity a {\bf degenerate CR singularity}. Some answers to the 
question just asked are known when $F$ is a globally-defined, proper branched covering
$F:\cplx\to\cplx$; refer to \cite{o'farrellPreskenis:uap2f89}. It would, however, be useful to know
what happens when $F$ is defined locally, i.e. to deduce whether $\srf={\rm graph}(F)$ is locally 
polynomially convex --- as in, for instance, the Forstneri{\v{c}}-Stout paper --- at the 
degenerate CR singularity $(0,F(0))$ just from local information about $F$. In a somewhat different
direction, Wiegerinck's paper \cite{wiegerinck:lpchdcrs95} studies {\em the failure of polynomial 
convexity} based on local conditions on $F$. In a recent paper \cite{bharali:sdCRslpc05}, the 
surface expressed locally at a degenerate CR singularity as
\begin{align}
\srf \ : \quad w \ &= \ \sum_{\alpha+\beta=k}C_{\alpha,\beta}z^\alpha{\zbar}^\beta+G(z) \notag\\
             &\equiv \ C_{k,0}z^k + C_{0,k}{\zbar}^{k} + \Sigma(z) + G(z), \label{E:greq}
\end{align}
where $k>2$ and $G$ is a smooth function satisfying $G(z)=o(|z|^{k})$ as $z\to 0$, was considered. 
A rough statement of one of the results in \cite{bharali:sdCRslpc05} is that given $C_{0,k}\neq 0$, if 
\begin{equation}\label{E:early}
\sup_{|\zt|=1}\frac{|\Sigma(\zt)|}{|\zt|^k} \ < \ 
	|C_{0,k}|\min\left\{\frac{\pi}{2k}, \frac{1}{2} \ \right\},
\end{equation}
and if $\Sigma(z)$ does not fluctuate too greatly, then $\srf$ is locally polynomially 
convex at $(0,0)$. The analytical condition \eqref{E:early} essentially says 
that --- writing $\srf={\rm graph}(F)$ locally --- if the $\zbar^k$ term is in some sense
the dominant term among all the leading-order terms in the Taylor expansion of $F$ around $z=0$, then
$\srf$ is polynomially convex in a small neighbourhood of $(0,0)$. One can make the following 
observations about the result under discussion:
\begin{itemize}
\item[a)] One might ask what can be deduced if some term other than the $\zbar^k$ term is the dominant
term among all the leading-order terms in the Taylor expansion of $F$ around $z=0$. From that 
perspective, the hypothesis discussed above is somewhat restrictive.
\smallskip

\item[b)] A careful examination of the proof of \cite[Theorem 1]{bharali:sdCRslpc05} reveals that
even under the restrictive hypothesis that $\zbar^k$ be the dominant term among all the leading-order
terms of $F$, \cite[Theorem 1]{bharali:sdCRslpc05} can be strengthened. 
\end{itemize}  
Item (a) above is an issue requiring care because, for example, if $F$ were a homogeneous polynomial 
in $z$ and $\zbar$ of degree $k$ , and a term of the form $z^m\zbar^{k-m}$, with $m\geq(k-m)>0$, 
significantly dominated all other terms, then $\srf:={\rm graph}(F)$ would {\em not} be locally
polynomially convex at $(0,0)$. This issue is resolved in Theorem \ref{T:pcvx} below. In the 
process, this theorem furnishes a considerably broader sufficient condition for local polynomial
convexity. However, we need to define some further notation. If $\phi$ is a complex-valued function
that is of class $\smoo^1$ on an open region $\OM\subset\cplx$, we define
\[
\|\nabla\phi(\zt)\| \ := \left|\partl{\phi}{z}(\zt)\right| \ + \ 
		\left|\partl{\phi}{\zbar}(\zt)\right| \quad\forall\zt\in\OM.
\]
We can now state our next theorem.
\smallskip

\begin{theorem}\label{T:pcvx} Let $\srf$ be a smooth surface in $\CC$ that is described
near $(0,0)\in\srf$ by
\begin{align}
\srf \ : \quad w \ &= \ \sum_{j=0}^k C_jz^{k-j}{\zbar}^j+G(z) \notag\\
             &\equiv \ C_0z^k + \Sigma(z) + G(z), \label{E:eqn}
\end{align}
where $k>2$, $G$ is a function of class $\smoo^1$ around $z=0$, and $G=O(|z|^{k+1})$.
Define the set $I(\srf):=\{j\in\mathbb{N}:k/2<j\leq k \ \text{and $C_j\neq 0$}\}$ and,
for each $j\in I(\srf)$, define
$\tau_j(z):=(\Sigma(z)-C_jz^{k-j}{\zbar}^j)/C_jz^{k-j}{\zbar}^j \ \forall z\neq 0$.
Suppose $I(\srf)\neq\emptyset$ and that there exists an integer $M\in I(\srf)$ such that
\begin{equation}\label{E:size}
\sup_{|\zt|=1}|\tau_M(\zt)| \ < \ \frac{\tan\{\pi/(2M-k)\}}{1+\tan\{\pi/(2M-k)\}},
\end{equation}
and
\begin{equation}\label{E:deriv}
k\sup_{|\zt|=1}|\tau_M(\zt)| + 
\left\{1-\sup_{|\zt|=1}|\tau_M(\zt)|\right\}^{-1}\sup_{|\zt|=1}|\zt| \ \|\nabla\tau_M(\zt)\|
\ < \ 2M-k.
\end{equation}
Then, there exists a small constant $\eps>0$ such that $\srf\cap\cyl{\eps}$ is polynomially
convex. Furthermore, calling the function on the right-hand side of \eqref{E:eqn} $\Phi$,
we have $\smoo(\ball{0}{\eps})=[z,\Phi]_{\ball{0}{\eps}}$.
\end{theorem}

\noindent{{\bf Remark.} It might seem on comparison that, owing to the condition \eqref{E:deriv}, 
for the case $M=k$, the above theorem is weaker than \cite[Theorem 1]{bharali:sdCRslpc05},
where the case $M=k$ has been treated. However, the bound \eqref{E:early} used in that result
is so stringent that, in fact, \eqref{E:deriv} is automatically implied. The conditions
of the above theorem are thus more permissive.}
\smallskip

\section{Some remarks on our proof techniques}\label{S:lemmas}

The primary purpose of this section is to state two brief lemmas that we shall need to 
prove the two theorems stated in \S\ref{S:intro}. In doing so, we shall also make a few
comments about the broad steps involved in the proofs of our theorems. 
\smallskip

To begin with, we state a lemma due to Bishop. This lemma may be found in the first half
of the proof of Theorem 4 in the paper \cite{bishop:mbfa59}. Before stating it, we note
that for the remainder of this paper, $\Leb$ will denote the Lebesgue measure on $\mathbb{C}$.
Bishop's lemma is vital to the proof of Theorem \ref{T:uaprx}, and is as follows:
\smallskip

\begin{lemma}[Bishop]\label{L:Bish} Let $\Dsc$ be a closed disc in $\cplx$. For any
measure $\mu\in\smoo(\Dsc)^\star$, define
\[
h_{\mu}(\zt) \ := \ \int_{\Dsc}\frac{d\mu(z)}{z-\zt}.
\]
Then, $|h_{\mu}|<\infty$ $\Leb$-a.e. in $\cplx$. If $h_{\mu}=0$ $\Leb$-a.e. in $\cplx$,
then $\mu=0$.
\end{lemma}

With $F$ as in Theorem \ref{T:uaprx}, we note that by definition
$[z,F]_{\Dsc}$ is a closed $\cplx$-linear subspace of $\smoo(\Dsc)$. Therefore, there exist
complex measures $\mu\in\smoo(\Dsc)^\star$ representing those continuous linear functionals on 
$\smoo(\Dsc)$ that annihilate $[z,F]_{\Dsc}$. Such measures will be called {\bf annihilating
measures}. The strategy behind the proof of Theorem \ref{T:uaprx} is to show that any 
annihilating measure $\mu$ is the zero measure. To achieve this, we fix an annihilating
measure $\mu$ in the proof of Theorem \ref{T:uaprx}, and carry out the following two steps:
\begin{itemize}
\item For each fixed $\zt$ lying off a certain exceptional set $\exep\varsubsetneq\Dsc$ with
$\Leb(\exep)=0$, we use the condition \eqref{E:sector} to construct a sequence 
$\{f_n\}_{n\in\mathbb{N}}\subset[z,F]_{\Dsc}$ such that $f_n(z)\longrightarrow 1/(z-\zt)$
$\mu$-a.e., and such that these functions are dominated by a function in 
$\mathbb{L}^1(d|\mu|,\Dsc)$.
\smallskip

\item Next, we apply the dominated convergence theorem to $\{f_n\}_{\nat}$ to show that 
$h_\mu(\zt)=0$ for each $\zt\notin\exep$, from which --- in view of Bishop's 
lemma --- Theorem \ref{T:uaprx} follows.
\end{itemize}
Details of these steps are presented in \S\ref{S:uaprx}.
\smallskip

The second lemma that plays an important role in this paper is Kallin's lemma. This is a
device that is used to determine when a union of polynomially convex sets is polynomially
convex. We state a certain form of Kallin's lemma that we shall use in \S\ref{S:pcvx}; the reader
is referred to \cite{kallin:fpcs66} for Kallin's original result.
\smallskip

\begin{lemma}[Kallin]\label{L:K-S} Suppose $X_1$ and $X_2$ are compact subsets of $\cplx^n$
such that $\poly(X_j)=\smoo(X_j)$, $j=1,2$. Let $\phi:\cplx^n\to\cplx$ be a holomorphic
polynomial such that $\phi(X_j)\subset W_j$, $j=1,2$, where $W_1$ and $W_2$ are polynomially
convex compact sets in $\cplx$ and $W_1\cap W_2=\{0\}$. Assume that
$\phi^{-1}\{0\}\cap(X_1\cup X_2)=X_1\cap X_2$. Then $\poly(X_1\cup X_2)=\smoo(X_1\cup X_2)$.
\end{lemma}

The above version of Kallin's lemma is presented within the proof of Theorem IV in 
\cite{forstnericStout:ncpcs91}. The symbol $\poly(K)$ denotes the uniform closure on
$K$ (where $K$ is compact) of the polynomials in $z$. The above lemma is useful in 
the context of the technique used in our proof of Theorem \ref{T:pcvx}. This proof will
essentially consist of the following three steps:
\begin{itemize}
\item{\em Step I: } We find a proper mapping $\Psi:\CC\to\CC$ and a $\delta>0$ such that
$\Psi^{-1}(\srf\cap\cyl{\delta})$ is a union of bordered surfaces $\presrf_1(\delta),\dots,
\presrf_{2M-k}(\delta)$, each of which enjoys certain special properties.
\smallskip

\item{\em Step II: } Each surface $\presrf_j(\delta)$ is the graph of a function $F_j$,
$j=1,\dots,2M-k$, over $\ball{0}{\delta}$. We show, using the properties $F_j$ possesses, that 
Theorem \ref{T:uaprx} is applicable to each $[z,F_j]_{\ball{0}{\delta}}$, $j=1,\dots,2M-k$.
\smallskip

\item{\em Step III: } Finally, we use Kallin's lemma to find an $\eps\in(0,\delta)$ such that 
$\poly(\presrf_1(\eps)\cup\dots\cup\presrf_{2M-k}(\eps))=
\smoo(\presrf_1(\eps)\cup\dots\cup\presrf_{2M-k}(\eps))$. Since $\Psi$ is a proper
covering, the last conclusion can be translated, using standard arguments, into the
conclusion of Theorem \ref{T:pcvx}.  
\end{itemize}
Details of the above argument are presented in \S\ref{S:pcvx}.
\smallskip

\section{The proof of Theorem \ref{T:uaprx}}\label{S:uaprx}
\smallskip

The following lemma is central to proving Theorem \ref{T:uaprx}.
\smallskip

\begin{lemma}\label{L:key} Let $\unal$ be a uniform algebra on a closed disc $\Dsc\Subset\cplx$ that
contains the function $z$. Fix $\zt\in\Dsc$. Assume that there is a function $W\subset\unal$ with 
the property that the set
\[
S_\zt \ := \ \{(z-\zt)W(z):z\in\Dsc\setminus(\{\zt\}\cup W^{-1}\{0\}) \ \}
\]
is contained in some open sector $S(0;I)$ with vertex at $0\in\cplx$. Then, there exists a
sequence of functions $\{f_n\}_{n\in\nat}\subset\unal$ such that
\[
\lim_{n\to\infty}f_n(z) \ = \ \frac{1}{z-\zt} 
			\quad\forall z\in\Dsc\setminus(\{\zt\}\cup W^{-1}\{0\}),
\]
and
\[
|f_n(z)| \ \leq \ \frac{4}{|z-\zt|} 
			\quad\forall z\in\Dsc\setminus(\{\zt\}\cup W^{-1}\{0\}), \
			\text{and $\forall n\in\nat$}.
\]
\end{lemma}
\noindent{{\em Proof.} We can find a $\phi\in[-\pi,\pi)$ and an integer $\nu\in\{1,2\}$ such 
that
\[
\er[(e^{i\phi}w)^{1/\nu}] \ > \ 0 \quad\forall w\in S(0;I),
\]
where the $\nu^{{\rm th}}$-root above is the appropriate branch of the $\nu^{{\rm th}}$-root that
is analytic on $\cplx\setminus(-\infty,0]$ and achieves the above inequality. Notice that it either 
suffices to choose $\nu=1$, or that $\nu=2$ necessarily, depending on whether the interior angle of
$S(0;I)$ is at most $\pi$, or is strictly greater than $\pi$. Define the holomorphic functions:
\[
P_n(w) \ := \ \left\{1-\frac{1}{[1+(e^{i\phi}w)^{1/\nu}]^n}\right\}\frac{1}{(e^{i\phi}w)^{1/\nu}}
\quad\forall w\in S(0;I), \ \forall n\in\nat.
\]}
Notice that $P_n$ extends to a continuous function on $\sect{0}{I}$. Therefore, defining
\[
Q_n(w) \ := \ \begin{cases}
		e^{i\phi}P_n(w)^\nu \ = \ 
		\left\{1-\dfrac{1}{[1+(e^{i\phi}w)^{1/\nu}]^n}\right\}^\nu\dfrac{1}{w} \ , &
		\text{if $w\in\sect{0}{I}\setminus\{0\}$},\\
		{} & {} \\
		e^{i\phi}n, & \text{if $w=0$},
		\end{cases}
\]
we conclude that $Q_n\in\hol(S(0;I))\cap\smoo(\sect{0}{I})$ for each $n\in\nat$.
\smallskip

Since $|1+(e^{i\phi}w)^{1/\nu}|>1 \ \forall w\in S(0;I)$, we have
\begin{equation}\label{E:lim}
\lim_{n\to\infty} Q_n(w) \ = \ 1/w \quad\forall w\in S(0;I)
\end{equation}
and, for the same reason
\begin{equation}\label{E:bd}
|Q_n(w)| \ \leq \ \left|1+\frac{1}{|1+(e^{i\phi}w)^{1/\nu}|^n}\right|^\nu\frac{1}{|w|} \ \leq
	\ \frac{4}{|w|} \quad\forall w\in S(0;I), \ \forall n\in\nat.
\end{equation}
Recall that:
\begin{itemize}
\item Since $W$ is continuous, there exists an $R>0$ such that 
$\overline{S}_\zt\subseteq\sect{0}{I}\cap\ball{0}{R}$;
\smallskip

\item  $Q_n\in\hol[S(0;I)\cap D(0;R)]\cap\smoo[\sect{0}{I}\cap\ball{0}{R}]$ for each 
$n\in\nat$.
\end{itemize}
Since $\sect{0}{I}\cap\ball{0}{R}$ is simply connected, by Mergelyan's theorem each
$Q_n$ is uniformly approximable on the set $\sect{0}{I}\cap\ball{0}{R}$ by polynomials in $z$.
Therefore, if we define $q_n(z):=Q_n[(z-\zt)W(z)]$, then $q_n\in\unal$ for each $n\in\nat$.
Now define
\[
f_n(z) \ := \ W(z)q_n(z) \quad\forall z\in\Dsc.
\]
Clearly $\{f_n\}_{n\in\nat}\subset\unal$. Observe that
\[
\lim_{n\to\infty}f_n(z) \ = \ W(z)\left[ \ \lim_{n\to\infty}Q_n((z-\zt)W(z)) \ \right]
		\ = \ \frac{1}{z-\zt}
		\quad\forall z\in\Dsc\setminus(\{\zt\}\cup W^{-1}\{0\}),
\]
which follows from \eqref{E:lim}, since $(z-\zt)W(z)\in S(0;I) \ 
\forall z\in\Dsc\setminus(\{\zt\}\cup W^{-1}\{0\}).$ For the same reason, \eqref{E:bd}
implies that
\[
|f_n(z)| \ \leq \ |W(z)|\times\frac{4}{|(z-\zt)W(z)|} \ = \ \frac{4}{|z-\zt|}
\quad\forall z\in\Dsc\setminus(\{\zt\}\cup W^{-1}\{0\}).\qed
\]
\smallskip

We now have all the tools necessary to provide the
\smallskip

\noindent{{\bf The proof of Theorem \ref{T:uaprx}.} Let $\mu\in\smoo(\Dsc)^\star$ be a measure that
annihilates $[z,F]_{\Dsc}$ (see \S\ref{S:lemmas} for a definition). In view of Bishop's lemma,
i.e. Lemma \ref{L:Bish}, we need to show that for the chosen annihilating measure $\mu$,
$h_{\mu}=0$ $\Leb$-a.e. So, we first consider $\zt\notin\Dsc$. Then, there exists a 
sequence of polynomials $\{p_n\}_{n\in\nat}$ that approximates the function 
$z\longmapsto(z-\zt)^{-1}$ uniformly on $\Dsc$. Clearly, $\{p_n\}_{n\in\nat}\subset[z,F]_{\Dsc}$.
Thus, owing to uniform convergence:}
\begin{equation}\label{E:van1}
0 \ = \ \lim_{n\to\infty}\int_{\Dsc}p_n \ d\mu \ = \ \int_{\Dsc}\frac{d\mu(z)}{z-\zt}\qquad
(\zt\notin\Dsc).
\end{equation}

Now define $A:=\{a\in\Dsc\setminus E:\mu(\{a\})\neq 0\}$. Since $\mu$ is a finite, regular
measure, $A$ is countable. Hence the set
\[
\widetilde{E} \ := \ \bigcup_{a\in A}F^{-1}\{F(a)\}
\]
is, by hypothesis, a countable union of countable sets. Since $z\longmapsto|z|^{-1}$ is locally
integrable with respect to the Lebesgue measure, and since $\mu$ is a finite measure supported 
in $\Dsc$, if we define 
\[
H_\mu(\zt) \ := \ \int_{\Dsc}\frac{d|\mu|(z)}{|z-\zt|} \ ,
\]
then $H_\mu<\infty$ $\Leb$-a.e. on $\Dsc$. Let 
$E^*=\{\zt\in\Dsc:H_\mu(\zt)=\infty\}$. Define
\[
\exep \ := \ E\bigcup\widetilde{E}\bigcup E^*.
\]
By the discussion just concluded, $\Leb(\exep)=0$. Now pick a $\zt\in\Dsc\setminus\exep$. Define
$W(z):=(z-\zt)\{F(z)-F(\zt)\}$. By \eqref{E:sector}, Lemma \ref{L:key} is applicable
with this choice of $W$ and with $\unal=[z,F]_{\Dsc}$. Thus, there exists a sequence of functions
$\{f_n\}_{n\in\nat}\subset\unal$ such that:
\begin{equation}\label{E:limbis}
\lim_{n\to\infty}f_n(z) \ = \frac{1}{z-\zt} 
			\quad\forall z\in\Dsc\setminus F^{-1}\{F(\zt)\},
\end{equation}
and
\begin{equation}\label{E:bdbis}
|f_n(z)| \ \leq \ \frac{4}{|z-\zt|} 
			\quad\forall z\in\Dsc\setminus F^{-1}\{F(\zt)\}, \ 
			\text{and $\forall n\in\nat$}.
\end{equation}
Note that since $\zt\notin(E\cup\widetilde{E})$, $\mu(F^{-1}\{F(\zt)\})=0$. Thus
\[
\eqref{E:limbis} \ \Longrightarrow \ f_n(z)\longrightarrow \frac{1}{z-\zt} 
					\quad\text{$\mu$-a.e.}
\]
Furthermore, as $\zt\notin E^*$, $H_\mu(\zt)<\infty$. Thus, in the present context:
\[
\eqref{E:bdbis} \ \Longrightarrow \ \text{The functions $f_n$ are dominated $\mu$-a.e.
					by a function in $\mathbb{L}^1(d|\mu|;\Dsc)$.}
\]
In view of the last two assertions, we may apply the dominated convergence theorem
as follows:
\begin{equation}\label{E:van2}
0 \ = \ \lim_{n\to\infty}\int_{\Dsc}f_n \ d\mu \ = \ \int_{\Dsc}\frac{d\mu(z)}{z-\zt}\qquad
(\zt\in\Dsc\setminus\exep).
\end{equation}

From \eqref{E:van1} and \eqref{E:van2}, we conclude that 
$h_\mu(\zt)=0 \ \forall\zt\in\cplx\setminus\exep$. This means that $h_\mu=0$ $\Leb$-a.e., whence 
$\mu=0$. Since this is true for any annihilating measure, $[z,F]_{\Dsc}=\smoo(\Dsc)$.
\begin{flushright} \qed \end{flushright}
\smallskip

\section{The proof of Theorem \ref{T:pcvx}}\label{S:pcvx} 
\smallskip

Before proceeding with the proof of our second theorem, we clarify two pieces of notation that
we shall use in the following proof. The expression $\phi\in\smoo^1(\ball{0}{\eps})$ will signify
that $\phi$ has continuous first-order derivatives at all points in some neighbourhood
of the closed disc $\ball{0}{\eps}$. On the other hand, the expression
$\phi\in\smoo^1(\ball{0}{\eps}^*)$ will mean that $\phi$ has continuous first-order derivatives
at all points in some neighbourhood of $\ball{0}{\eps}$ {\em except at $0\in\mathbb{C}$}.
\smallskip

Next, we define a couple of concepts that will be used in the proof below. Firstly, 
if $K$ is a compact subset of $\cplx^n$, the {\bf polynomially convex hull of $\boldsymbol{K}$}, 
written $\widehat{K}$, is defined
by
\[
\widehat{K} \ := \ \{\zt\in\cplx^n \ | \ |P(\zt)|<\sup_K|P|, \ \text{for every holomorphic
                        polynomial $P$}\}.
\]
Given a uniform algebra $\mathscr{A}$, the {\bf maximal ideal space of $\boldsymbol{\mathscr{A}}$}
is the space of all unit-norm algebra-homomorphisms of $\mathscr{A}$ to $\cplx$, viewed as a subset 
of the dual space $\mathscr{A}^\star$ with the $\text{weak}^\star$ topology 
(it is a standard fact that every
complex homomorphism of $\mathscr{A}$ is in fact continuous). Recall that for a compact subset $K$,
the maximal ideal space of $\smoo(K)$ is homeomorphically identified with $K$. We will need this 
fact in the following proof.
\smallskip

Having established these preliminaries, we are in a position to give
\smallskip

\noindent{{\bf The proof of Theorem \ref{T:pcvx}.} We begin by introducing a new system of
global holomorphic coordinates $(\zee,\wee)$ defined by
\begin{align}
\zee \ &:= \ z, \notag\\
\wee \ &:= \ w-C_0z^k. \notag
\end{align}
Relative to these new coordinates, $\srf$ is expressed as
\[
\srf \ : \quad \wee \ = \ \Sigma(\zee)+G(\zee).
\]
For simplicity of notation, we shall denote the new coordinates by $(z,w)$, and work with the 
following presentation of $\srf$:
\begin{equation}\label{E:clean}
\srf \ : \quad w \ = \ \Sigma(z)+G(z),
\end{equation}
where the meanings of $\Sigma$ and $G$ remain unchanged from those in \eqref{E:eqn}. Moreover,
the reader may check that neither of the hypotheses \eqref{E:size} or \eqref{E:deriv} are
affected by this change of coordinate. Let us refer to the right-hand side of \eqref{E:clean}
by $\fee(z)$.
\smallskip

Let $I(\srf)$ be as in Theorem \ref{T:pcvx}, and let $M\in I(\srf)$ be such that the
associated $\tau_M$ satisfies the conditions \eqref{E:size} and \eqref{E:deriv}. Define
$\pow:=2M-k$. Observe that $\pow>0$. Define the map $\Psi:\CC\to\CC$ by $\Psi(z,w):=(z,w^\pow)$. This
is a proper, holomorphic, $\pow$-to-$1$ covering map. We now present the first step of this proof.
\smallskip

\noindent{{\bf Step I:} {\em We show that there exists a small constant $\delta>0$ such that
$\Psi^{-1}(\srf\cap\cyl{\delta})=\bigcup_{j=1}^\pow\presrf_j(\delta)$, where, 
for $0<r\leq\delta$, $\presrf_j(r)$ represent graphs of the form
\begin{equation}\label{E:pregraphs}
\presrf_j(r) \ := \ \{(z,w):w=C_*\omega_j\{|z|^{k/\pow}e^{-i\tht}+f(z)+R(z)\}, \ 
				|z|\leq r\}, \quad j=1,\dots,\pow,
\end{equation}
and where
\begin{itemize}
\item We write $z:=|z|e^{i\tht}$;
\item$C_*:=|C_M|^{1/\pow}e^{i{\rm Arg}(C_M)/\pow}$, and $\omega_j=e^{2\pi i(j-1)/\pow}$, i.e.
a $\pow^{th}$-root of unity;
\item $f\in\smoo^{1}(\ball{0}{\delta})$ if $M\neq k$, but $f\in\smoo^{1}(\ball{0}{\delta}^*)$ 
if $M=k$; and
\item $R\in\smoo^{1}(\ball{0}{\delta})$ and $R(z)=O(|z|^{1+(k/\pow)})$.
\end{itemize}}}

\noindent{To see this, we first note that if $\pow=1$ then $f(z)=(\Sigma(z)-C_Mz^{k-M}\zbar^M)/C_M$
and $R(z)=G(z)/C_M$, and the stated properties of $f$ and $G$ are obvious from our hypotheses.}
Therefore, we may focus on the case $\pow\neq 1$. In this situation, we first write
\begin{equation}\label{E:fee}
\fee(z) \ := \ C_Mz^{k-M}\zbar^M\left\{1+\tau_M(z)+\frac{G(z)}{C_Mz^{k-M}\zbar^M}\right\}, \quad
z\neq 0.
\end{equation}
Observe that
\begin{itemize}
\item Owing to the estimate \eqref{E:size}, $|\tau_M(z)|<1$; and
\item $G(z)/C_Mz^{k-M}\zbar^M\longrightarrow 0$ uniformly as $z\to 0$.
\end{itemize}
For these reasons, we can find a small constant $\delta>0$ such that
\[
\left|\tau_M(z)+\frac{G(z)}{C_Mz^{k-M}\zbar^M}\right| \ < \ 1 \quad\forall z:0<|z|\leq\delta.
\]
Given this fact, $\fee(z)$ has $\pow$ distinct $\pow^{th}$-roots --- call them 
$F_j(z), \ j=1,\dots,\pow$ --- when $0<|z|\leq\delta$, which are obtained by applying the 
Binomial Theorem with exponent $1/\pow$ to the expression in \eqref{E:fee} that is enclosed
in brackets. This results in the expression
\begin{align}
F_j(z) \ &= \ C_*\omega_j\left\{|z|^{k/\pow}e^{-i\tht}+
\sum_{\nu=1}^\infty\alpha_\nu|z|^{k/\pow}e^{-i\tht}\tau_M(z)^\nu+R(z)\right\}\notag\\
&\equiv \ C_*\omega_j\{|z|^{k/\pow}e^{-i\tht}+f(z)+R(z)\}, \quad\forall z:|z|\leq\delta.\notag
\end{align}
Note that in the above expression, the quantity $|z|^{k/\pow}e^{-i\tht}\tau_M(z)$ is interpreted
as
\[
|z|^{k/\pow}e^{-i\tht}\tau_M(z) \ = \ \begin{cases}
					|z|^{k/\pow}e^{-i\tht}\tau_M(z), &\text{if $0<|z|\leq\delta$},\\
					0, &\text{if $z=0$},
					\end{cases}
\]
and the $\alpha_\nu$'s are the coefficients occurring in the Taylor expansion of $(1+x)^{1/\pow}$
around $x=0$.
\smallskip

Note that since the series expansion that produces $F_j(z)$ converges absolutely and uniformly 
for the specified range of $z$, rearrangement is permissable; and it is by rearrangement that
the quantities $f(z)$ and $R(z)$ are constructed. Furthermore, term-by-term differentiation is
possible. But note that 
\[
|z|^{k/\pow}e^{-i\tht}\tau_M(z) \ \text{is} \ \begin{cases}
				\text{continuously differentiable at $z=0$ if $k>\pow$,}\\
				\text{not differentiable at $z=0$ if $k=\pow$.}
				\end{cases}
\]
The degrees of regularity claimed for $f$ and $R$ readily follow from the last two statements.
This concludes the first step of our proof.
\smallskip

\noindent{{\bf Step II:} {\em We show that there exists an $\eps\in(0,\delta]$ such that
$\poly(\presrf_j(\eps))=\smoo(\presrf_j(\eps)), \ j=1,\dots,\pow$.}

\noindent{Note that each $\presrf_j(\delta)$ is a complex linear image of the set
$\{(z,w):w=|z|^{k/\pow}e^{-i\tht}+f(z)+R(z), \ |z|\leq\delta\}$. Therefore, to accomplish
this step, it suffices to find an $\eps>0$ such that 
$[z,F_0]_{\ball{0}{\eps}}=\smoo(\ball{0}{\eps})$, where $F_0(z):=|z|^{k/\pow}e^{-i\tht}+f(z)+R(z)$.
We will use Theorem \ref{T:uaprx} to accomplish this.}
\smallskip

We shall first need a few computations. Writing $z=|z|e^{i\tht}=re^{i\tht}$, recall that
\[
\partl{{}}{z} \ = \ \frac{e^{-i\tht}}{2}\left\{\partl{{}}{r}-\frac{i}{r}\partl{{}}{\tht}\right\},
\qquad
\partl{{}}{\zbar} \ = \ \frac{e^{i\tht}}{2}\left\{\partl{{}}{r}+\frac{i}{r}\partl{{}}{\tht}\right\}.
\]
Therefore, we have
\begin{align}\label{E:prtlMain}
\partl{{}}{z}( \ |z|^{k/\pow}e^{-i\tht}) \ &= \ \frac{e^{-2i\tht}}{2}
\left(\frac{k}{\pow}-1\right)|z|^{(k/\pow)-1}, \\
\partl{{}}{\zbar}( \ |z|^{k/\pow}e^{-i\tht}) \ &= \ \frac{1}{2}
\left(\frac{k}{\pow}+1\right)|z|^{(k/\pow)-1}. \notag
\end{align}
Another computation that we will find useful is the following estimate for the quantity
\[
\left|\partl{f}{z}(\zt)\right| +
\left|\partl{f}{\zbar}(\zt)\right|, \quad\zt\neq 0. 
\]
Using the expressions \eqref{E:prtlMain}, we compute:
\begin{align}
\left|\partl{f}{z}(\zt)\right| &+
\left|\partl{f}{\zbar}(\zt)\right| \notag \\
&\leq \ 
\frac{1}{2}\left(\frac{k}{\pow}-1\right)|\zt|^{(k/\pow)-1}
\left|\sum_{\nu=1}^\infty\alpha_\nu\tau_M(\zt)^\nu\right| \ + \ 
|\zt|^{k/\pow}\sum_{\nu=1}^\infty\nu|\alpha_\nu||\tau_M(\zt)|^{\nu-1}
\left|\partl{\tau_M}{z}(\zt)\right| \notag \\
&\quad + \ \frac{1}{2}\left(\frac{k}{\pow}+1\right)|\zt|^{(k/\pow)-1}
\left|\sum_{\nu=1}^\infty\alpha_\nu\tau_M(\zt)^\nu\right| \ + \
|\zt|^{k/\pow}\sum_{\nu=1}^\infty\nu|\alpha_\nu||\tau_M(\zt)|^{\nu-1}
\left|\partl{\tau_M}{\zbar}(\zt)\right| \notag \\
&\leq \frac{k}{\pow} \ |\zt|^{(k/\pow)-1}\sum_{\nu=1}^\infty|\alpha_\nu||\tau_M(\zt)|^\nu 
\ + \ |\zt|^{k/\pow}\|\nabla\tau_M(\zt)\|\sum_{\nu=1}^\infty\nu|\alpha_\nu||\tau_M(\zt)|^{\nu-1}\notag \\
&= \ \frac{k}{\pow} \ |\zt|^{(k/\pow)-1}\left[1-(1-|\tau_M(\zt)| \ )^{1/\pow}\right] \notag \\
&\hskip0.5in + \ 
|\zt|^{k/\pow}\|\nabla\tau_M(\zt)\|
\left.\frac{d}{dx}\left[1-(1-x)^{1/\pow}\right]\right|_{x=|\tau_M(\zt)|}\quad (0<|\zt|\leq\delta)\notag \\
&\leq \ \frac{k}{\pow} \ |\zt|^{(k/\pow)-1}|\tau_M(\zt)| \ + \ 
\frac{|\zt|^{k/\pow}}{\pow}\|\nabla\tau_M(\zt)\| \ (1-|\tau_M(\zt)| \ )^{-1}
\quad (0<|\zt|\leq\delta)\label{E:gradEst}
\end{align}
The last line of the above estimate follows from the fact that as $|\tau_M(\zt)|<1$ for the relevant
range of $\zt$, we have the inequalities $(1-|\tau_M(\zt)| \ )<(1-|\tau_M(\zt)| \ )^{1/\pow}<1$.
\smallskip

In the language of Theorem \ref{T:uaprx}, define $E:=\{0\}$. The above computations will allow us
to determine the distribution of values of the quantity
\[
(z-\zt)\{F_0(z)-F_0(\zt)\} \quad\zt\notin E, \ z\notin F_0^{-1}\{F_0(\zt)\},
\]
for each $\zt\in\ball{0}{\eps}\setminus E$ and for all $z\in\ball{0}{\eps}\setminus F_0^{-1}\{F_0(\zt)\}$,
for some sufficiently small $\eps\in(0,\delta]$. For this purpose, we:
\begin{itemize}
\item {\em Fix} a $\zt:0<|\zt|\leq\delta$;
\item Define $A:=\sup_{|\zt|=1}|\tau_M(\zt)|$; and
\item Define
\[
B \ := \ \frac{k}{\pow} \ \sup_{|\zt|=1}|\tau_M(\zt)| +\frac{1}
{\pow}\left\{1-\sup_{|\zt|=1}|\tau_M(\zt)|\right\}^{-1}\sup_{|\zt|=1}|\zt| \ \|\nabla\tau_M(\zt)\|.
\]
\end{itemize}
Our task will have to be taken up under two different cases.
\smallskip

\noindent{{\bf{\em Case 1:}} {\em EITHER $k/2<M<k$, OR $M=k$ and the line 
joining $z$ to $\zt$ does not contain the origin in its interior}.}

\noindent{We will explain the reason behind this unusual division into cases in a moment. First,
however, we define the real-valued function $\psi_{\zt,z}:[0,1]\longrightarrow\mathbb{R}$ by
\[
\psi_{\zt,z}(t) \ := \ \er[\ (z-\zt)F_0(tz+(1-t)\zt) \ ].
\]
Note that when $k/2<M<k$, then for $|z|$ sufficiently small, the $\psi_{\zt,z}$ thus defined
would be of class $\smoo^1$ in a small neighbourhood of the unit-interval. In the present
setting, even when $M=k$, $\psi_{\zt,z}$ is of class $\smoo^1$ on the open unit-interval. Note
that since $f$ is not differentiable at the origin, {\em the preceding statements would not be
true if $M=k$ and the line joining $z$ to $\zt$ contained the origin in its interior.} This 
explains the necessity of the present division into cases. By this
discussion, we see that the Mean Value Theorem is applicable to $\psi_{\zt,z}$. Thus, there
exists a $t_*\in(0,1)$ such that
\begin{equation}\label{E:realPre}
\er[ \ (z-\zt)\{F_0(z)-F_0(\zt)\} \ ] \ = \ \psi_{\zt,z}(1)-\psi_{\zt,z}(0) \ 
= \ \psi_{\zt,z}^\prime(t_*).
\end{equation}
Define $\xi_*:=t_*z+(1-t_*)\zt$. It {\em may be possible} that $\xi_*=0$. In that case
\begin{equation}\label{E:real1}
\xi_*=0 \ \Longrightarrow \ \begin{cases}
				\er[ \ (z-\zt)\{F_0(z)-F_0(\zt)\} \ ] \ = \ 0, \\
				(z-\zt)\{F_0(z)-F_0(\zt)\} \ \neq \ 0.
				\end{cases}
\end{equation}
The second statement above follows from the fact that $z\notin F_0^{-1}\{F_0(\zt)\}$.
Now, we consider the case $\xi_*\neq 0$. We use the calculations \eqref{E:prtlMain} and
\eqref{E:gradEst} to compute
\begin{align}
\psi_{\zt,z}^\prime(t_*) \ &= \ \er\left[ \ (z-\zt)^2\partl{F_0}{z}(\xi_*)+
|z-\zt|^2\partl{F_0}{\zbar}(\xi_*)\right] \notag\\
&\geq \ \frac{1}{2}\left(\frac{k}{\pow}+1\right)|\xi_*|^{(k/\pow)-1}|z-\zt|^2
-\frac{1}{2}\left(\frac{k}{\pow}-1\right)|\xi_*|^{(k/\pow)-1}|z-\zt|^2 \notag\\
&\quad\qquad -\left\{ \ \left|\partl{f}{z}(\xi_*)\right| + \left|\partl{f}{\zbar}(\xi_*)\right| \ 
\right\}|z-\zt|^2 - O( \ |\xi_*|^{k/\pow}|z-\zt|^2) \notag \\
&\geq \ |\xi_*|^{(k/\pow)-1}|z-\zt|^2-|\xi_*|^{(k/\pow)-1}B|z-\zt|^2-O( \ |\xi_*|^{k/\pow}|z-\zt|^2)
\quad(\xi_*\neq 0)\label{E:real2}
\end{align}
The final inequality above follows from \eqref{E:gradEst} and from the definition of $B$.
Note that by the condition \eqref{E:deriv}, $B<1$. Thus, we can find a constant $\eps_1>0$ 
so small that
\begin{multline}
0<|\zt|\leq\eps_1, \ z\in\ball{0}{\eps_1}\setminus F_0^{-1}\{F_0(\zt)\} \ \text{and}
\ \xi_*\neq 0 \\
\Longrightarrow \  |\xi_*|^{(k/\pow)-1}|z-\zt|^2-|\xi_*|^{(k/\pow)-1}B|z-\zt|^2
-O( \ |\xi_*|^{k/\pow}|z-\zt|^2) \ > \ 0.\notag
\end{multline}
It is pertinent to note here that {\em even though} $\xi_*$ itself may vary somewhat
unpredictably as $(z,\zt)$ is varied, the manner in which $\xi_*$ enters the estimate
\eqref{E:real2} makes it possible to choose $\eps_1>0$ uniformly in the above estimate.
Combining this estimate with \eqref{E:realPre} and \eqref{E:real2}, we get
\begin{multline}\label{E:real3}
0<|\zt|\leq\eps_1, \ z\in\ball{0}{\eps_1}\setminus F_0^{-1}\{F_0(\zt)\} \ \text{and} 
\ \xi_*\neq 0 \\
\Longrightarrow \  \psi_{\zt,z}^\prime(t_*) \ = \ \er[ \ (z-\zt)\{F_0(z)-F_0(\zt)\} \ ] \ > \ 0.
\end{multline}
From \eqref{E:real1} and \eqref{E:real3}, we conclude that
\begin{align}
\text{\em Under the conditions of Case 1, there exists a constant $\eps_1$}&
\text{$\in(0,\delta]$ such that}\notag \\
\er[ \ (z-\zt)\{F_0(z)-F_0(\zt)\} \ ] \ \in \ \{\xi\in\cplx:\er(\xi)\geq 0, \ \xi\neq 0\}, \ 
&z\in\ball{0}{\eps_1} \setminus F_0^{-1}\{F_0(\zt)\},\notag \\
&\forall\zt\in\ball{0}{\eps_1}\setminus E.\label{E:real4}
\end{align}
\smallskip

\noindent{{\bf{\em Case 2:}} {\em $M=k$ and the line joining $z$ to $\zt$ contains the origin 
in its interior}.}

\noindent{The analysis of this case is an obvious variation of the method used in Step 1; hence
we shall be brief. We define two functions $\psi_z$ and $\psi_\zt : [0,1]\to\mathbb{R}$ by
\begin{align}
\psi_{z}(t) \ &:= \ \er[\ (z-\zt)F_0(tz) \ ], \notag \\
\psi_{\zt}(t) \ &:= \ \er[\ (z-\zt)F_0((1-t)\zt) \ ]. \notag
\end{align}
Both functions are continuous on $[0,1]$ and differentiable in $(0,1)$. Thus, by the Mean Value
Theorem, there exist $\tau_1, \ \tau_2\in (0,1)$ such that
\[
\psi_{z}(1)-\psi_{z}(0) \ = \ \psi_{z}^\prime(\tau_1),\qquad 
\psi_{\zt}(1)-\psi_{\zt}(0) \ = \ \psi_{\zt}^\prime(\tau_2).
\]
Note that
\begin{equation}\label{E:preReal}
\er[ \ (z-\zt)\{F_0(z)-F_0(\zt)\} \ ] \ = \ [\psi_{z}(1)-\psi_{z}(0)] \ + \ 
[\psi_{\zt}(1)-\psi_{\zt}(0)] \ = \ \psi_{z}^\prime(\tau_1) + \psi_{\zt}^\prime(\tau_2).
\end{equation}
The next observation is vital to our estimates of $\psi_{z}^\prime(\tau_1)$ and
$\psi_{\zt}^\prime(\tau_2)$. Since $z$ and $\zt$ lie on a line through the origin,
and on {\em opposite sides of the origin},
\begin{equation}\label{E:vital}
(z-\zt)\zbar \ = \ |z||z-\zt|, \qquad -\overline{\zt}(z-\zt) \ = \ |\zt||z-\zt|.
\end{equation}
Now, define $\xi_1:=\tau_1z$, and $\xi_2:=\tau_2\zt$. We emulate the calculations 
leading up to \eqref{E:real2} above, in conjunction with the relations in \eqref{E:vital},
to get
\begin{align}
\psi_{z}^\prime(\tau_1) \ &= \ \er\left[z(z-\zt)\partl{F_0}{z}(\xi_1)+
\zbar(z-\zt)\partl{F_0}{\zbar}(\xi_1)\right] \notag\\
&\geq \ |z||z-\zt|-\left\{ \ \left|\partl{f}{z}(\xi_1)\right| + \left|\partl{f}{\zbar}(\xi_1)\right| \
\right\}|z||z-\zt| - O( \ |\xi_1|^{k/\pow}|z||z-\zt| \ ) \notag \\
&\geq \ |z||z-\zt|-B|z||z-\zt|-O( \ |\xi_1||z||z-\zt| \ )
\label{E:real5}
\end{align}
An analogous calculation gives
\begin{equation}\label{E:real6}
\psi_{\zt}^\prime(\tau_2) \ \geq \ 
|\zt||z-\zt|-B|\zt||z-\zt|-O( \ |\xi_2||\zt||z-\zt| \ ).
\end{equation}
Arguing exactly as in Case 1, we can infer from \eqref{E:preReal},, \eqref{E:real5} and
\eqref{E:real6} that:
\begin{align}
\text{\em Under the conditions of Case 2, there exists a constant $\eps_2$}&
\text{$\in(0,\delta]$ such that}\notag \\
\er[ \ (z-\zt)\{F_0(z)-F_0(\zt)\} \ ] \ \in \ \{\xi\in\cplx:\er(\xi)>0\}, \ 
&z\in\ball{0}{\eps_2} \setminus F_0^{-1}\{F_0(\zt)\},\notag \\
&\forall\zt\in\ball{0}{\eps_2}\setminus E.\label{E:real7}
\end{align}
\smallskip

Note that \eqref{E:real4} and \eqref{E:real7} verify condition \eqref{E:sector} of
Theorem \ref{T:uaprx} for $F_0$ with $E=\{0\}$ and $\Dsc=\ball{0}{\min(\eps_1,\eps_2)}$.
It remains to examine the cardinality of $F_0^{-1}\{F_0(\zt)\}$ for $\zt\neq 0$. Notice
that by the very same considerations as in the calculation \eqref{E:gradEst}, we obtain
\begin{align}
&\left|\partl{F_0}{\zbar}(\zt)\right|-
\left|\partl{F_0}{z}(\zt)\right| \notag \\
&\qquad\geq \ \frac{1}{2}\left(\frac{k}{\pow}+1\right)|\zt|^{(k/\pow)-1} - 
\frac{1}{2}\left(\frac{k}{\pow}-1\right)|\zt|^{(k/\pow)-1} \notag \\
&\qquad\qquad\qquad-\left\{ \ \left|\partl{f}{z}(\zt)\right| 
+\left|\partl{f}{\zbar}(\zt)\right| \ \right\}
-O( \ |\zt|^{k/\pow}) \notag \\
&\qquad\geq \ |\zt|^{(k/\pow)-1} - |\zt|^{(k/\pow)-1}\left\{\frac{k}{\pow} \ |\tau_M(\zt)| \ 
+ \ \frac{|\zt|}{\pow}\|\nabla\tau_M(\zt)\| \ (1-|\tau_M(\zt)| \ )^{-1}\right\}
-O( \ |\zt|^{k/\pow}) \notag \\
&\qquad\geq \ |\zt|^{(k/\pow)-1} - B|\zt|^{(k/\pow)-1} - O( \ |\zt|^{k/\pow}) \notag
\end{align}
Once again, as $B<1$, we can find a constant $\eps_3$ that is so small that 
\begin{equation}\label{E:Jac}
\left|\partl{F_0}{\zbar}(\zt)\right|-\left|\partl{F_0}{z}(\zt)\right| \ \geq \ 
(1-B)|\zt|^{(k/\pow)-1} - O( \ |\zt|^{k/\pow}) \ > \ 0 \ \forall\zt:0<|\zt|\leq\eps_3.
\end{equation}
We can view $F_0$ as a mapping of the disc $\ball{0}{\delta}$ in $\mathbb{R}^2$ into
$\mathbb{R}^2$. Therefore, we can define the real Jacobian of this mapping (except at
$\zt=0$ when $M=k$), which we denote by ${\rm Jac}_{\mathbb{R}}(F_0)$. The inequality
\eqref{E:Jac} tells us that
\[
{\rm Jac}_{\mathbb{R}}(F_0)(\zt) \ = \ 
\left|\partl{F_0}{z}(\zt)\right|^2-\left|\partl{F_0}{\zbar}(\zt)\right|^2 \ < \ 0 \quad 
\forall\zt:0<|\zt|\leq\eps_3.
\]
By the Inverse Function Theorem, we conclude from the above statement that for each
$\zt\in D(0,\eps_3)\setminus\{0\}$, $F_0^{-1}\{F_0(\zt)\}\cap D(0;\eps_3)$ is a 
discrete set. Thus, if we define
\[
\eps \ := \ \min\{\eps_1, \ \eps_2, \ \eps_3/2\},
\]
then, both the hypotheses of Theorem \ref{T:uaprx} are satisfied for $F_0$ with 
$E=\{0\}$ and $\Dsc=\ball{0}{\eps}$. Hence, $[z,F_0]_{\ball{0}{\eps}}=\smoo(\ball{0}{\eps})$.
Now, note that
\[
\poly(\presrf_j(\eps)) =  [z,C_*\omega_jF_0]_{\ball{0}{\eps}} =  
[z,F_0]_{\ball{0}{\eps}} = \smoo(\ball{0}{\eps}) = \smoo(\presrf_j(\eps)), \quad
j=1,\dots,\pow.
\]
Hence, the second step of our proof is accomplished.
\smallskip

The following step, as we shall see, completes the proof of Theorem \ref{T:pcvx}.
\smallskip

\noindent{{\bf Step III:} {\em We apply Kallin's Lemma to the conclusions of Step II to 
prove the desired result.}

\noindent{Consider the polynomial $p(z,w)=zw/C_*$. For any $(z,w)\in \presrf_1(\eps)$,
\begin{align}
\er\{p(z,w)\} \ &= \ |z|^{(k/\pow)+1} + 
\er\left\{ \  |z|^{(k/\pow)+1}\sum_{\nu=1}^\infty\alpha_\nu\tau_M(z)^\nu + zR(z)\right\} \notag \\
&\geq \ |z|^{(k/\pow)+1} - |z|^{(k/\pow)+1}\sum_{\nu=1}^\infty|\alpha_\nu||\tau_M(z)|^\nu
- O( \ |z|^{(k/\pow)+2}) \notag\\
&\geq \ |z|^{(k/\pow)+1} - |z|^{(k/\pow)+1}\left[1-(1-|\tau_M(z)| \ )^{1/\pow}\right]
- O( \ |z|^{(k/\pow)+2}) \notag\\
&\geq \ |z|^{(k/\pow)+1} - A|z|^{(k/\pow)+1} - O( \ |z|^{(k/\pow)+2}).\label{E:repart}
\end{align}
Similarly, for any $(z,w)\in\presrf_1(\eps)$, we estimate
\begin{align}
|\mi\{p(z,w)\}| \ &\leq \ |z|^{(k/\pow)+1}\left\{ \ \sum_{\nu=1}^\infty|\alpha_\nu||\tau_M(z)|^\nu +
|R(z)| \ \right\}.\notag \\
&\leq \ A|z|^{(k/\pow)+1} + O( \ |z|^{(k/\pow)+2}).\label{E:impart} 
\end{align}
Recall that $A:=\sup_{|\zt|=1}|\tau_M(\zt)|$. Let us fix a constant $C$ such that
\begin{equation}\label{E:angle}
A \ < \ C \ < \frac{\tan(\pi/\pow)}{1+\tan(\pi/\pow)}.
\end{equation}
Examining the expressions \eqref{E:repart} and \eqref{E:impart}, we see that we can, lowering the value
of $\eps>0$ if necessary, arrange for
\begin{align}
|\mi\{p(z,w)\}| \ &\leq \ C|z|^{(k/\pow)+1}, \notag \\
\er\{p(z,w)\} \ &\geq \ (1-C)|z|^{(k/\pow)+1}, \quad\forall (z,w)\in\presrf_1(\eps).\notag
\end{align}
Note that lowering the value of $\eps>0$ {\em does not} alter the conclusion of Step II above.
In view of the last inequalities:
\[
p(\presrf_1(\eps)) \ \varsubsetneq \ 
\left\{x+iy\in\cplx : |y|\leq \frac{C}{1-C}x, \ x\geq 0 \right\}.
\]
The above expression says that $p(\presrf_1(\eps))$ is a proper subset of the closed sector
$W_1$ that is centered on the positive $x$-axis, and has an vertex-angle of 
\[
2\arctan\left(\frac{C}{1-C}\right) \ < \ 2\arctan\{\tan(\pi/\pow)\} \ = \ 2\pi/\pow.
\]
The above inequality is a consequence of the condition \eqref{E:angle} on $C$.
Note that by construction, when $j\neq 1$, $p(\presrf_j(\eps))$ is a proper subset of the 
closed sector $W_j$, which is simply a copy of $W_1$ rotated by $(2\pi(j-1)/\pow)$, $j=2,\dots,\pow$.
We have shown so far that:
\begin{itemize}
\item For each $\presrf_j(\eps)$, $\poly(\presrf_j(\eps))=\smoo(\presrf_j(\eps))$, $j=1,\dots,\pow$;
\item $p(\presrf_j(\eps))\varsubsetneq W_j, \ j=1,\dots,\pow$;
\item $W_\mu\cap W_\nu = \{0\} \ \forall \mu\neq\nu$, because the vertex-angle of each
$W_j$ is less than $2\pi/\pow$; and
\item $p^{-1}\{0\}\cap\left\{ \ \bigcup_{j=1}^\pow \presrf_j(\eps) \ \right\} = \{(0,0)\}$.
\end{itemize}
The above facts allow us to apply Kallin's Lemma repeatedly to show that
\begin{equation}\label{E:sff_jcon}
\poly\left(\cup_{j=1}^\pow\presrf_j(\eps)\right) \ = \ 
\smoo\left(\cup_{j=1}^\pow\presrf_j(\eps)\right).
\end{equation}
Now let $\psi\in\smoo(\cyl{\eps}\cap\srf)$. Define
$\widehat{\psi}:=\psi\circ\Psi : \Psi^{-1}(\cyl{\eps}\cap\srf)\to\cplx$. As
$\Psi^{-1}(\cyl{\eps}\cap\srf)=\bigcup_{j=1}^k\presrf_j(\eps)$,
$\widehat{\psi}\in\smoo\left(\bigcup_{j=1}^\pow\presrf_j(\eps)\right)$. We can paraphrase 
\eqref{E:sff_jcon} in the following way:
for each $\sma>0$, there exists a polynomial $g_\sma$ such that
\begin{equation}\label{E:presrf_japprox}
|\widehat{\psi}(z,e^{2\pi i(j-1)/\pow}w)-g_\sma(z,e^{2\pi i(j-1)/\pow}w)| \ < \ \sma
\quad \forall (z,w)\in \presrf_1(\eps),
\quad j=1,\dots,\pow.
\end{equation}
We define
\[
Q_\sma(z,w) \ :=  \ \frac{1}{\pow}\sum_{j=1}^\pow g_\sma(z,e^{2\pi i(j-1)/\pow}w).
\]
Notice that if $g_\sma(z,w)=\sum_{0\leq\mu+\nu\leq N}A_{\mu,\nu}z^\mu w^\nu$, then
$Q_\sma(z,w)$ has the form
\begin{align}
Q_\sma(z,w) \ &= \ \sum_{(\mu,\nu):\nu=\pow j}A_{\mu,\pow j}z^\mu w^{\pow j} \notag\\
        &\equiv \ P_\sma(z,w^\pow), \notag
\end{align}
where $P_\sma$ is itself a polynomial. Let us write $w=|w|e^{i\phi}$, $\phi\in[0,2\pi)$.
For $(z,w)\in\cyl{\eps}\cap\srf$, we compute
\begin{align}
|\psi(z,w)&-P_\sma(z,w)| \ = \ \left|\frac{1}{\pow}
\sum_{j=1}^\pow\widehat{\psi}(z,|w|^{1/\pow}e^{i(2\pi(j-1)+\phi)/\pow}) -
Q_\sma(z,|w|^{1/\pow}e^{i\phi/\pow})\right|\notag\\
&\leq \ \sum_{j=1}^\pow \ \frac{|\widehat{\psi}(z,|w|^{1/\pow}e^{i(2\pi(j-1)+\phi)/\pow}) -
g_\sma(z,|w|^{1/\pow}e^{i(2\pi (j-1)+\phi)/\pow})|}{\pow} \ < \ \pow\left(\frac{\sma}{\pow}\right).\notag
\end{align}
The last inequality follows from the estimate \eqref{E:presrf_japprox}. This
establishes that $\poly(\cyl{\eps}\cap\srf) = \smoo(\cyl{\eps}\cap\srf)$. Equivalently, we
have just established (recall that $\Phi$ denotes the right-hand side of \eqref{E:eqn}) that 
$[z,\Phi]_{\ball{0}{\eps}}=\smoo(\ball{0}{\eps})$.
\smallskip
        
We now only need to show that $\cyl{\eps_0}\cap\srf$ is polynomially convex. This follows from
general abstract considerations. For this purpose, given a compact $K\Subset \mathbb{C}^n$, we
define
\begin{align}
\widehat{K} \ &:= \ \text{the polynomially convex hull of $K$},\notag\\
\mathscr{P}(K;\mathbb{C}^n) \ &:= \ \text{the uniform algebra on $K$ generated by the class
                                        $\{f|_K:f\in\cplx[z_1,\dots,z_n] \ \}$},\notag\\
\mideal \ &:= \ \text{the maximal ideal space of the uniform algebra
$\mathscr{P}(K;\mathbb{C}^n)$}.\notag
\end{align}
It is well known that $\mideal=\widehat{K}$; see, for instance, \cite[Chap.6/\S29]{stout:tua71}.
Thus, in our situation, 
$\mathscr{M}[\mathscr{P}(\cyl{\eps}\cap\srf;\CC)]=\{(z,w):|z|\leq\widehat{\eps\}\cap}\srf$.
But since $\poly(\cyl{\eps}\cap\srf) = \smoo(\cyl{\eps}\cap\srf)$,
\[
\{(z,w):|z|\leq\widehat{\eps\}\cap}\srf \ = \ \mathscr{M}[\smoo(\cyl{\eps_0}\cap\srf)] \ = \
\cyl{\eps}\cap\srf.
\]
This concludes our proof.
\begin{flushright} \qed \end{flushright}

\end{document}